\documentclass[12pt,reqno,a4wide]{amsart}

\usepackage{tikz} 
\usepackage{calrsfs}
\usepackage{mathrsfs}

\usepackage{hyperref}

\usetikzlibrary{decorations.pathreplacing}
\usetikzlibrary{fit}

\usepackage{pgfplots}

\allowdisplaybreaks

\newtheorem{theorem}{Theorem}

\catcode`,\active

\catcode`\,12

\begin{document}
\bibliographystyle{plain}

%
%

	\title[A matrix with sums of Catalan  numbers]
	{A matrix with sums of Catalan  numbers---LU-decomposition and determinant }

	\author[H. Prodinger ]{Helmut Prodinger }
	\address{Department of Mathematics, University of Stellenbosch 7602, Stellenbosch, South Africa
	and
	Department of Mathematics and Mathematical Statistics,
	Umea University,
	907 36 Umea, 	 Sweden  }
	\email{hproding@sun.ac.za}

	
	\begin{abstract}
Following Benjamin et al., a matrix with entries being sums of two neighbouring Catalan numbers is considered. Its LU-decomposition is
given, by guessing the results and later prove it by computer algebra, with lots of human help. Specializing a parameter, the determinant
turns out to be a Fibonacci number with odd index, confirming earlier results, obtained back then by combinatorial methods.
	\end{abstract}
	

\maketitle


%
%

\section{Introduction}

Let $\mathscr{C}_n=\frac1{n+1}\binom{2n}{n}$ be the $n$-th Catalan number.
The $n\times n$ Matrix 
\begin{equation*}
\mathscr{M}=\left(\begin{matrix}
	\mathscr{C}_{t}+\mathscr{C}_{t+1}&\mathscr{C}_{t+1}+\mathscr{C}_{t+2}&\dots&\mathscr{C}_{t+n-1}+\mathscr{C}_{t+n}\\
	\mathscr{C}_{t+1}+\mathscr{C}_{t+2}&\mathscr{C}_{t+2}+\mathscr{C}_{t+3}&\dots&\mathscr{C}_{t+n}+\mathscr{C}_{t+n+1}\\
\vdots&\vdots&\ddots&\vdots\\
\mathscr{C}_{t+n-1}+\mathscr{C}_{t+n}&\mathscr{C}_{t+n}+\mathscr{C}_{t+n+1}&\dots&\mathscr{C}_{t+2n-2}+\mathscr{C}_{t+2n-1}\\
\end{matrix}\right)
\end{equation*}
is considered in  \cite{Yerger}; the determinant is considered by combinatorial means. The natural range of the parameters is $n\ge1$ and $t\ge0$.
There are many methods to compute determinants of combinatorial matrices, as expertly 
described in \cite{kratt1, kratt2}.

In this paper, we consider the LU-decomposition $LU=\mathscr{M}$, with a lower triangular matrix $L$ 
with 1's on the main diagonal, and an upper triangular matrix $U$. From this, the determinant
comes out as a corollary, by multiplying the elements in $U$'s main diagonal. We restrict our attention 
to the instance $t=0$, since the computations seem to become very messy in the more general setting. But at the same time, we
consider a more general matrix with an extra parameter $x$, viz.
 \begin{equation*}
 	\mathscr{M}=\left(\begin{matrix}
 		\mathscr{C}_{0}+x\mathscr{C}_{1}&\mathscr{C}_{1}+x\mathscr{C}_{2}&\dots&\mathscr{C}_{n-1}+x\mathscr{C}_{n}\\
 		\mathscr{C}_{1}+x\mathscr{C}_{2}&\mathscr{C}_{2}+x\mathscr{C}_{3}&\dots&\mathscr{C}_{n}+x\mathscr{C}_{n+1}\\
 		\vdots&\vdots&\ddots&\vdots\\
 		\mathscr{C}_{n-1}+x\mathscr{C}_{n}&\mathscr{C}_{n}+x\mathscr{C}_{n+1}&\dots&\mathscr{C}_{2n-2}+x\mathscr{C}_{2n-1}
 	\end{matrix}\right).
 \end{equation*}
Not only do we get more general results in this way, but it is actually easier to guess the explicit forms of $L$ and $U$ with an extra parameter involved.

Here are the results that we found by computer experiments, which we consider to be the main contributions of this paper: 

\begin{theorem}
	For $k,i\ge1$, set
\begin{align*}
	F(k,i)&=\frac1{i(2i-1)}\binom{2i}{i-k}\sum_{0\le r\le k}\frac{1}{2k-r}\binom{2k-r}{r}\Big(ri+2ik^2-ik-2rk^2+2k^3-k^2\Big)x^r
\end{align*}
and
\begin{equation*}
	g(k)=\sum_{0\le r\le k}\binom{2k-r}{r}x^r=F(k,k).
\end{equation*}
Then
\begin{equation*}
L[i,k]=\frac{F(k,i)}{g(k)}\quad\text{and}\quad U[k,j]=\frac{F(k,j)}{g(k-1)}.
\end{equation*}
\end{theorem} 

In the next section, first the expressions for $F(i,j)$ and $g(k)$ will be simplified, and then it will be proved that these two matrices are indeed the
LU-decomposition of $\mathscr{M}$. Note that only \emph{one} function $F(k,i)$ is used to represent both, $L[i,k]$ and $U[k,j]$. This shows in particular the symmetry
related to $i\leftrightarrow j$.

\section{Simplification and proof}

In many instances where Catalan numbers are involved, it is beneficial to work with an auxiliary variable:
\begin{equation*}
x=\frac{-u}{(1+u)^2}\quad\text{and}\quad u=\frac{-1-2x+\sqrt{1+4x}}{2x}.
\end{equation*}
Then 
\begin{equation*}
g(k)=\frac{1-u^{2k+1}}{(1-u)(1+u)^{2k}}.
\end{equation*}
This is well within the reach of modern computer algebra (I use Maple). Further,
\begin{equation*}
	F(k,j)=	(1- u^{2k}	) \frac{ \binom{2 j}{	j-k
	}}{2j (2 j-1)} \frac{2k^2-j}{(1-u)(1+u)^{2k-1}}
	+(1+u^{2k}	) \frac{  \binom{2 j}{	j-k
		}k}{2j  (1+u)^{2k}}.
\end{equation*}
Maple is capable to simplify $F(k,j)$, but the version given here, which is pleasant, was obtained with help from Carsten Schneider and his software~\cite{carsten}.
Of course, once this version is known, Maple can confirm that it is equivalent to its own simplification. Note that $F(k,k)=g(k)$, and the L-matrix has indeed
1's on the main diagonal.

What is nice to note is that $L[i,k]=0$ for $i<k$ and $U[k,j]=0$ for $k>j$ automatically, thanks to the properties of binomial coefficients: a binomial coefficient $\binom nm$
with integers $n,m$ such that $n\ge0$ and $m<0$ is equal to zero.

Now we want to evaluate the $(i,j)$ entry of the matrix $L\cdot U$:
\begin{equation*}
\sum_{k\ge1}L[i,k]U[k,j].
\end{equation*}
Maple cannot evaluate this sum without help:
\begin{equation*}
\frac{F(k,i)F(k,j)}{g(k)g(k-1)}=\frac{\text{expression}}{(1-u^{2k+1})(1-u^{2k-1})}
\end{equation*}
What helps here is partial fraction decomposition:
\begin{equation*}
	\frac{F(k,i)F(k,j)}{g(k)g(k-1)}=\text{expression}_1+\frac{\text{expression}_2}{(1-u^{2k+1})}
	+\frac{\text{expression}_3}{(1-u^{2k-1})}.
\end{equation*}
In the second term the change of index $k\to k-1$ makes things better, so that Maple can  compute the sum over $k$;
however, a correction term needs to be taken in:
\begin{equation*}
	\sum_{k=1}^j\frac{F(k,i)F(k,j)}{g(k)g(k-1)}=\sum_{k=1}^j\frac{\text{expression}_4}{(1-u^{2k-1})}-\frac{\text{expression}_2}{(1-u^{2k+1})}\Big|_{k=0}.
\end{equation*}
All the expressions are long and can be created with a computer. The sum can now be computed, and, switching back to the $x$-world,
simplifies (again with a lot of human help, e. g., to simplify expressions in which the Gamma-functions appears)  the last sum to
\begin{equation*}
\mathscr{C}_{i+j-2}+x\mathscr{C}_{i+j-1},
\end{equation*}
as it should. For our simplification, we still used the variable $u$ in~\cite{HPmaple}. However, for small $x$ and $u$, the connection between the two variables is bijective.

All the details can be checked in the maple worksheet~\cite{HPmaple}. Perhaps a quick comment how the partial fraction decomposition is working is the essential formula
\begin{equation*}
\frac1{g(k)g(k-1)}=(1-u)(1+u)^{4k-3}\Big[\frac{1}{1-u^{2k-1}}-\frac{u^2}{1-u^{2k+1}}\Big].
\end{equation*}

\section{The determinant}

The values in the main diagonal are given by
\begin{equation*}
U[k,k]=\frac{g(k)}{g(k-1)}.
\end{equation*}
Consequently
\begin{equation*}
	\prod_{k=1}^nU[k,k]=\frac{g(n)}{g(0)}=g(n).
\end{equation*}
Setting $x=1$, as in \cite{Yerger}, means $u=-\frac{3+\sqrt5}{2}=-\alpha^2$, with $\alpha=\frac{1+\sqrt5}{2}$ being the
golden ratio. We also need $\beta=\frac{1-\sqrt5}{2}$. After some straightforward simplifications, this can be rewritten in terms of Fibonacci numbers:
\begin{equation*}
g(n)=\frac{1+\alpha^{4n+2}}{(1-\alpha^2)^{2n}(1+\alpha^2)}=\frac{1+\alpha^{4n+2}}{\alpha^{2n}\sqrt 5\alpha}
=\frac{\alpha^{2n+1}-\beta^{2n+1}}{\sqrt5}=F_{2n+1}.
\end{equation*}

\bibliographystyle{plain}

\end{document}